\documentclass[12pt, a4paper,reqno]{amsart}

\usepackage[T1]{fontenc}
\usepackage[utf8]{inputenc}
\usepackage{xcolor}

\definecolor{bleu_sombre}{rgb}{0,0,0.6}  \definecolor{rouge_sombre}{rgb}{0.8,0,0}\definecolor{vert_sombre}{rgb}{0,0.6,0}
\usepackage[plainpages=false,colorlinks,linkcolor=bleu_sombre,
citecolor=rouge_sombre,urlcolor=vert_sombre,breaklinks]{hyperref}

\usepackage[english]{babel}
\usepackage{geometry}

\usepackage{amsmath,amssymb,amsthm,graphicx,amsfonts,enumerate,stmaryrd, mathrsfs,bbm}

\theoremstyle{plain}
\newtheorem{theorem}{{Theorem}}[section]
\newtheorem*{theorem*}{{Theorem}}
\newtheorem{proposition}[theorem]{Proposition}

\newtheorem*{proposition*}{Proposition}

\newtheorem*{corollary*}{Corollary}
\newtheorem{lemma}[theorem]{Lemma}

\newtheorem*{lemma*}{Lemma}

\theoremstyle{definition}

\newtheorem*{definition*}{Definition}

\theoremstyle{remark}
\newtheorem{remark}[theorem]{Remark}
\newtheorem{notation}[theorem]{Notation}

\makeatletter

\@addtoreset{equation}{section}
\makeatother

\renewcommand{\leq}{\leqslant}	\renewcommand{\geq}{\geqslant}

\newcommand{\R}{\mathbb{R}}

\renewcommand{\Re}{\mathrm{Re}\,}
\renewcommand{\Im}{\mathrm{Im}\,}

\title[Purely magnetic tunneling]{Purely magnetic tunneling \\ between radial magnetic wells}

\author[S. Fournais]{S\o ren Fournais}
\address[S. Fournais]{Department of Mathematics, University of Copenhagen, Universitetsparken 5, DK-2100 Copenhagen \O, Denmark}
\email{fournais@math.ku.dk}

\author[L. Morin]{Léo Morin}
\address[L. Morin]{Department of Mathematics, University of Copenhagen, Universitetsparken 5, DK-2100 Copenhagen \O, Denmark}
\email{leo.morin@math.au.dk }

\author[N. Raymond]{Nicolas Raymond}
\address[N. Raymond]{Univ Angers, CNRS, LAREMA, Institut Universitaire de France, SFR MATHSTIC, F-49000 Angers, France}
\email{nicolas.raymond@univ-angers.fr}

\begin{document}
		\begin{abstract}
		This article is devoted to the semiclassical spectral analysis of the magnetic Laplacian in two dimensions. Assuming that the magnetic field is positive and has two symmetric radial wells, we establish an accurate tunneling formula, that is a one-term estimate of the spectral gap between the lowest two eigenvalues. This gap is exponentially small when the semiclassical parameter goes to zero, but positive.
	\end{abstract}
	\maketitle

\section{Introduction}

\subsection{Motivation}
This article deals with the spectral analysis of the purely magnetic Laplacian in two dimensions, especially with the effect of a symmetry of the magnetic field on the spectrum. More precisely, we are interested in this effect in the semiclassical limit, which is equivalent to the large magnetic field limit. In such a regime, much has been done in the last two decades to describe the bound states of the magnetic Schrödinger operators, see the books \cite{FH10} and \cite{Raymond} (for earlier references see the classical texts \cite{CFKS,AHS} and references therein). In \cite{Raymond}, it is underlined that a motivation is to understand tunneling effects induced by magnetic fields and to extend the Helffer-Sjöstrand theory to magnetic fields. This theory was originally developed in the context of purely electric Schrödinger operators, see \cite{HS84, HS85b}. Helffer and Sjöstrand also added a (rather small) magnetic field in \cite{HS87} and were able to establish a tunneling formula. This problem was reconsidered recently in \cite{FSW22}, where the authors were able to go beyond the smallness assumptions on the magnetic field, but without establishing an accurate expansion of the spectral gap. This motivated another work \cite{HK22} where the authors prove various estimates of the spectral gap by means of a rather lengthy analysis. In fact, even though the context of the present article is different, its analysis will cast a new light on \cite{FSW22, HK22} and provide the reader with a concise strategy very close to the original spirit developed in \cite{HS84, HS85b}.

It turns out that the adaptation of the Helffer-Sjöstrand theory to the \emph{purely} magnetic situation was more complex than expected. Even establishing the semiclassical expansions of the eigenvalues in cases without symmetries was already a challenge (see, for instance, \cite{RVN15} where Birkhoff normal forms are used). A key point in the study of the tunneling effect is to have an accurate description of the eigenfunctions. Such a description was given rather recently in \cite{BHR16} (in a multiscale context) and in \cite{BR20, GBNRVN21} in the case of a magnetic well\footnote{By this, we mean that the magnetic field has a unique minimum, which is non-degenerate.} in two dimensions. To prove that these WKB constructions are good approximations in general is still a widely open problem (see  \cite{GBNRVN21} which gives a first result in this direction). Concerning the tunneling estimates in the case of pure magnetic fields in two dimensions, we have to mention \cite{BHR22} where the Neumann Laplacian with constant magnetic field is considered in a symmetric domain and the spectral gap accurately described by means of a microlocal dimensional reduction to the boundary. This work inspired similar results such as \cite{HKS23}, but also \cite{AA23} in the case of magnetic fields vanishing along curves and \cite{FHK22} in the case of discontinuous magnetic fields.

The present article tackles a situation with symmetric magnetic wells. In this case, we know that the eigenfunctions are localized near the bottoms of the wells (see \cite{GBRVN21}), and not along curves as in \cite{BHR22} (what invites \emph{a priori} to use another strategy). Until now, except for the WKB constructions of \cite{BR20, GBNRVN21}, nothing was known in the direction of tunneling results for double-well magnetic fields.

\subsection{The purely magnetic tunneling result}
Let us now describe the framework of this article and our main result. Choose lengths $a,L>0$ with $a<\frac{L}{2}$ and magnetic field strengths $0<b_0 <b_1$.
Define the two points (with distance $L$)
\[ x_\ell=\left(-\frac{L}{2},0\right), \qquad x_r=\left(\frac{L}{2},0\right).\]
Let the left magnetic field $B_{\ell} \in C^{\infty}(\R^2,\R)$, have minimal value $b_0>0$ attained uniquely and non-degenerately at $x_{\ell}$ and be constant equal to $b_1 > b_0$ on $\{|x-x_{\ell}|>a\}$. 
Furthermore, we assume that $B_{\ell}(x) \in [b_0,b_1]$ for all $x$.
We will define the right magnetic field $B_{r}$ by symmetry and
assume the wells to be radial and symmetric. More precisely,
\begin{equation}
B_\ell \text{ is a smooth function of } |x-x_\ell|^2, \text{ and }
B_r(x_1,x_2) = B_\ell(-x_1,x_2).
\end{equation}

We consider a double well magnetic field $B \in C^{\infty}(\R^2,\R)$, with minimal value $b_0 >0$ attained (non-degenerately) at exactly the two points $x_\ell$ and $x_r$.
The magnetic field will be constant outside the disks of radius $a$, through the definition
\begin{equation}
	B(x) =
	\begin{cases}
		B_{\ell}(x), &\text{if } |x-x_\ell| \leq a,\\
		B_r(x), &\text{if } |x-x_r| \leq a,\\
		b_1, &\text{otherwise.}
	\end{cases}
\end{equation}

Two important quantities of the problem are the flux difference
\begin{equation}\label{def.M0}
M := \frac{1}{4\pi} \int_{\R^2} \big( B(x) -b_1 \big) \mathrm d x =\frac{1}{2\pi} \int_{\R^2} \big( B_{\ell}(x) -b_1 \big) \mathrm d x \leq 0,
\end{equation}
and the relative flux (through the disc of radius $\frac{L}{2}$)
\begin{equation}\label{eq.N}
0<N:=\frac{8|M|}{b_1 L^2} < \frac{a^2}{(L/2)^2} < 1\,.
\end{equation}

We consider a smooth vector potential $\mathbf{A}=(0,A_2)$ associated with $B$, and the magnetic Laplacian $\mathscr{L}_h=(-ih\nabla-\mathbf{A})^2$. We denote by $\lambda_1(h) \leq \lambda_2(h)$ the two smallest eigenvalues of $\mathscr{L}_h$. In this article we prove the following estimate on the spectral gap.

\begin{theorem}\label{thm.main}
Assume that $L> (2+\sqrt 6) a$.
Then, there exists a positive constant $c$  such that the following holds. We have
\[\lambda_2(h)-\lambda_1(h)\underset{h\to 0}{\sim}c e^{-S/h} h^{\frac{1}{2} + \frac{b_0}{2b_1}}\,,\]
where
\begin{equation}\label{eq.S0delta0}
S= \int_{0}^{L/2} \frac{1}{\pi r^2} \left(\int_{D(x_\ell,r)} B(y)\mathrm{d}y\right) r \mathrm{d} r +I\,,
\end{equation}
and
\begin{equation}\label{eq.I}
I=\frac{b_1 L^2}{4}\left(\frac{N-1}{2}+\sqrt{1-N}-N\ln\Big(1+\sqrt{1-N}\Big)\right)>0\,.
\end{equation}
	\end{theorem}

As far as the authors know, Theorem \ref{thm.main} is the only known result of its kind in the case of pure magnetic wells. Note that:
\begin{enumerate}
\item The first integral term in $S$ is an Agmon distance between the two wells, see Proposition \ref{prop.WKB} below. Contrary to the case of electric wells, we have here an additional positive term $I$ inside the exponential. This term is a purely magnetic effect, and vanishes when $N=1$.
\item The variations of the field also appears in the power $b_0 / b_1$ of $h$ in the prefactor.
\item The constant $c$ can be computed explicitly in terms of $B$, however it has no simple interpretation, see \eqref{eq:c}.
\item The assumption on $L$ can be slightly lightened, see inequality \eqref{eq.inegalite2}. This condition is here to ensure that the remainder term in Proposition \ref{prop.w} is small enough.
\item An alternative way of writing $S$, which could also be of interest is the following formula,
$$
S-I = \frac{-1}{\pi} \int_{\R^2} B(x) \mathbbm{1}_{\{|x-x_{\ell}| \leq \frac{L}{2}\}} \ln\Big(\frac{2|x-x_{\ell}|}{L}\Big)\,\mathrm{d}x.
$$
\end{enumerate}

\subsection{Organization and strategy}
The paper is organized as follows. In Section~\ref{sec.2}, we consider the one-well operator (on the left) called $\mathscr{L}_{h,\ell}$ and we recall the asymptotic expansions of its eigenvalues $(\mu_n(h))_{n\geq 1}$ obtained by Helffer and Kordyukov in \cite{HK11} (or in \cite{RVN15} with a more geometric point of view), see Proposition~\ref{prop.HK}. Then, we describe the groundstate eigenfunction $\phi_\ell$ by recalling the results of \cite{GBNRVN21}. This eigenfunction is radial and can be approximated by a WKB Ansatz, see Proposition \ref{prop.WKB}. This accurate description of the groundstate of the one-well operator is known to be a corner stone in the analysis of the tunneling effect. 

In Section \ref{sec.3}, we recall a standard (and rough) localization result: the eigenfunctions of the double-well operator are exponentially localized near $x_\ell$ and $x_r$. This allows to decouple the wells and we get the first naive estimate $\lambda_2(h)-\lambda_1(h)=\mathscr{O}(h^\infty)$. Then, we use $\phi_\ell$ and $\phi_r$ (up to a change of gauge) as quasimodes for the double-well operator $\mathscr{L}_h$. This gives an upper bound, see Proposition \ref{prop.goodub}. 

Section \ref{sec.4} is devoted to the interaction term $w_h$ governing the behavior of the spectral gap $\lambda_2(h)-\lambda_1(h)$, see Proposition \ref{prop.w} and \ref{prop.wh}. The integral formula \eqref{eq.w0} for $w_h$ is similar to that obtained in the electric case in \cite{HS84}. This formula could be used in the context of \cite{FSW22, HK22} and it could simplify this analysis. Note that, contrary to the situation in \cite{FSW22, HK22}, the interaction phase $\theta(0,x_2)$ is not linear and the non-linear term is directly related to the fact that the magnetic field is not constant.

 In Section \ref{sec.5}, we establish the one-term asymptotics of $w_h$, see Proposition \ref{prop.finalw}. To estimate the integral giving \eqref{eq.w0}, we use the following observation (done in \cite{FSW22} with a constant magnetic field): the eigenfunction $\phi_\ell$ has an explicit expression in the interaction zone (this comes from the fact that the magnetic field is radial and constant between the wells), see Lemma \ref{lem.explicit}. Combining this with the WKB approximation of Proposition \ref{prop.WKB}, we get an explicit integral to estimate, see \eqref{eq.WWW} and Lemma \ref{lem.Ch}. In this integral (with three variables $s_1$, $s_2$ and $y$), the phase is complex. That is why we look for the complex critical points in the variable $y$ (the vertical variable), see Sections \ref{sec.phasey} and \ref{sec.crit}. Then, we can use a complex translation of the integral in $y$ and apply the Laplace method, see Section \ref{sec.w1w2}. This reduces the analysis to an integral in two variables \eqref{eq.reducedint}, which can be estimated by the usual Laplace method, see the description of the real phase $F(s)$ in Section \ref{sec.crit2}.

\section{The one-well operator and its groundstate}\label{sec.2}
We consider the vector potential
\begin{equation}\label{eq.Al}
\mathbf{A}_\ell(x)=\int_0^1 B_\ell(x_\ell +t(x-x_\ell))t(x-x_\ell)^{\perp}\mathrm{d}t\,,
\end{equation}
which is associated with $B_\ell$ and we consider the corresponding magnetic Laplacian \[\mathscr{L}_{h,\ell}=(-ih\nabla-\mathbf{A}_\ell)^2\,.\]
We denote by $(\mu_n(h))_{n\geq 1}$ the non-decreasing sequence of its eigenvalues. The following proposition is well-known.

\begin{proposition}[\cite{HK11} \& \cite{GBNRVN21}]\label{prop.HK}
	Let $n\geq 1$. For $h$ small enough, $\mu_n(h)$ belongs to the discrete spectrum of $\mathscr{L}_{h,\ell}$ and
	\[\mu_n(h)=b_0h+((2nd_0+d_1)h^2+o(h^2)\,,\]
	where
	\[d_0=\frac{\sqrt{\det H}}{b_0}\,,\quad d_1=\frac{(\mathrm{Tr} H^{\frac12})^2}{2b_0}\,,\quad H=\frac12\mathrm{Hess} B_{\ell}(x_{\ell})\,.\]
	\end{proposition}
In the following we will denote $\mu_h = \mu_1(h)$. Let us now describe the groundstate $\phi_\ell$ of $\mathscr{L}_{h,\ell}$ and its WKB approximation. It is already known that $\phi_\ell$ is radial in the variable $r=|x-x_\ell |$ and does not vanish if $h$ is small enough, see \cite{GBNRVN21}. In particular, it solves an ODE (see \cite[Section 4]{GBNRVN21}),
\begin{equation}\label{eq.radial.schrodinger}
\left(-h^2r^{-1}\partial_rr\partial_r+\left(\frac{\alpha(r)}{r}\right)^2\right)\phi_\ell=\mu_h\phi_\ell\,,\quad \alpha(r)=\int_0^{\frac{r^2}{2}}\beta(u) \mathrm{d} u\,,
\end{equation}
where the smooth function $\beta : [0,+\infty)\to\R_{>0}$, the magnetic profile, is defined by
\begin{equation}\label{def.beta}
B_\ell(x) = \beta \Big( \frac{|x-x_\ell|^2}{2} \Big).
\end{equation} Our assumptions imply that $\beta'(0)>0$ and that $0$ is the unique minimum of $\beta$. Equation \eqref{eq.radial.schrodinger} can be interpreted as the eigenvalue equation of a radial \emph{electric} Schrödinger operator. This observation, combined with the WKB construction in \cite[Theorem 2.5]{GBNRVN21} (see also \cite{BR20}), implies that $\phi_\ell$ is approximated by a WKB Ansatz (see for instance \cite{Hel88,HS84, HS85b}).

\begin{proposition}\label{prop.WKB}
Let
\[\Phi_\ell(x)=\frac{1}{2}\int_{0}^{|x-x_\ell|^2/2}\frac{1}{v}\left(\int_0^{v} \beta(u)\mathrm{d}u\right)\mathrm{d}v\,.\]	
There exists a sequence of smooth radial (with center at $x_\ell$) functions $(a_{j})_{j\geq 0}$ on $\R^2$, with $a_{0}>0$ such that the following holds. For any given $p\in\mathbb{N}$, letting $a_{h,\ell}(x)=\sum_{j=0}^p h^j a_j(x)$, and considering
\[\phi^{\mathrm{WKB}}_\ell(x)=h^{-\frac12}a_{h,\ell}(x)e^{-\Phi_\ell(x)/h}\,,\]
we have
\[e^{\Phi_\ell(x)/h}\left(\mathscr{L}_{h,\ell}-\mu_h\right)\phi_\ell^{\mathrm{WKB}}=\mathscr{O}(h^{p+1})\,.\]
Moreover, for all $\epsilon\in(0,1)$, there exist $C, h_0>0$ such that, for all $h\in(0,h_0)$,
\begin{equation}\label{eq.Agmon}
\|(-ih\nabla-\mathbf{A}_\ell)(e^{(1-\epsilon)\Phi_\ell}\phi_\ell)\|\leq Ch\|\phi_\ell\|\,.
\end{equation}
We also have the approximation
\[e^{\Phi_\ell(x)/h}(\phi_{\ell}(x)-\phi_\ell^{\mathrm{WKB}}(x))=\mathscr{O}(h^{p+1})\,,\]
locally uniformly.
	\end{proposition}

\begin{remark}
We can rewrite $\Phi_\ell$ in terms of $B_\ell$. In particular we have
\[ \Phi_\ell(0) = \frac 1 2 \int_{0}^{L/2} \frac{1}{\pi r^2} \left(\int_{D(x_\ell,r)} B(y)\mathrm{d}y\right) r \mathrm{d} r, \]
which is the quantity appearing in the coefficient $S$ from Theorem \ref{thm.main}. It is the flux of an averaged version of $B$.
\end{remark}

\begin{remark}
By symmetry, we have a similar WKB approximation in the right well. Considering the vector potential
\[\mathbf{A}_r(x)=\int_0^1 B_r(x_r +t(x-x_r))t(x-x_r)^{\perp}\mathrm{d}t\,,\]
we denote by $\mathscr{L}_{h,r}$ the corresponding magnetic Laplacian.	Its positive normalized groundstate is denoted by $\phi_r$.	
\end{remark}

\section{First estimates on the tunneling effect}\label{sec.3}

\subsection{A naive localization estimate}

The following propositions are already known, see, for instance, the discussion in \cite[Section 1.4]{GBRVN21}.
\begin{proposition}
Let $C_0>0$. There exist $\alpha,C, h_0>0$ such that for all $h\in(0,h_0)$ and all eigenpair $(\lambda,\psi)$ with $\lambda\leq b_0h+C_0h^2$	we have
\[\int_{\R^2}e^{2\alpha \min(|x-x_r|,|x-x_\ell|)/h^{1/4}}|\psi|^2\mathrm{d}x\leq C\|\psi\|^2\,.\]	
\end{proposition}

This rather rough exponential decay of the eigenfunctions is enough to prove that the spectrum of the double-well operator is the superposition of the spectra of the one-well operators, modulo $\mathscr{O}(h^\infty)$. 
\begin{proposition}\label{prop.specgap}
	There exist $c,h_0>0$ such that, for all $h\in(0,h_0)$, we have
	\[|\lambda_2(h)-\lambda_1(h)|=\mathscr{O}(h^\infty)\,,\]
	and $\lambda_3(h)-\lambda_1(h)\geq ch^2$.
\end{proposition}

\subsection{A better upper bound of the splitting}

We now use the WKB approximation, Proposition \ref{prop.WKB}, to prove a better upper bound on the spectral gap. Let
\begin{equation}\label{eq.S1}
S_0=\Phi_\ell\left(\frac{L}{2}-a,0\right)=\frac{1}{2}\int_{0}^{(L-a)^2/2}\frac{1}{v}\left(\int_0^{v} \beta(u)\mathrm{d}u\right)\mathrm{d}v\,,
\end{equation}
which can be interpreted as an Agmon distance from the left well to the border of the right well.

\begin{notation}
In the following, we use the classical notation (from \cite{HS84, HS85b}) that $f(h)=\tilde{\mathscr{O}}(e^{-S_0/h})$ when, for all $\epsilon>0$, we have $f(h)={\mathscr{O}}(e^{-(S_0-\epsilon)/h})$.
\end{notation}

\begin{proposition}	\label{prop.goodub}
We have
	\[|\lambda_2(h)-\lambda_1(h)|=\tilde{\mathscr{O}}(e^{-S_0/h} )\,,\quad \lambda_1(h)=\mu(h)+\tilde{\mathscr{O}}(e^{-S_0/h} )\,.\]

	\end{proposition}
	
	\begin{proof}
Fix a small $\eta>0$. Let	 $\chi_\ell$ be a smooth cutoff function equal to $1$ on $\{x_1<L/2-a-\eta\}$ and $0$ on $\{x_1>L/2-a-\eta/2\}$. We let
\begin{equation}\label{eq.phitilde}
\tilde\phi_\ell(x)=\chi_\ell(x)\hat\phi_\ell(x)\,,\quad \hat\phi_\ell(x)=e^{i\sigma_\ell(x)}\phi_\ell(x)\,,
\end{equation}	
	where
	\[\sigma_\ell(x)=\int_{[0,x]}(\mathbf{A}-\mathbf{A}_\ell)\cdot\mathrm{d}s\]
is a gauge change from $\mathbf A_\ell$ to $\mathbf A$.
	Note that, on $\{x_1<L/2-a\}$
	\[\nabla\sigma_\ell=\mathbf{A}-\mathbf{A}_\ell\,.\]
With this choice we have,  on $\{x_1<L/2-a-\eta\}$,
\[(\mathscr{L}_h-\mu_h)\tilde\phi_\ell=0\,.\]
Moreover,  thanks to \eqref{eq.Agmon}, we have
\[\|(\mathscr{L}_h-\mu_h)\tilde\phi_\ell\|=\tilde{\mathscr{O}}(e^{-S_\eta/h})\,,\]
where
\begin{equation}\label{eq.Seta}
S_\eta = \Phi_\ell \Big( \frac L 2 - a - \eta, 0 \Big) = \Phi_r \Big( - \frac L 2 + a + \eta, 0 \Big).
\end{equation}	
In the same way, and using a similar notation for the right well, we have
\[\|(\mathscr{L}_h-\mu_h)\tilde\phi_r\|=\tilde{\mathscr{O}}(e^{-S_\eta/h})\,.\]	
Since $\tilde\phi_\ell$ and $\tilde\phi_r$ are linearly independent (for $h$ small enough), the spectral theorem implies that there are at least two eigenvalues of $\mathscr{L}_h$ close to $\mu_h$ modulo $\tilde{\mathscr{O}}(e^{-S_\eta/h})$. These are necessarily the first two eigenvalues since $\mu_h$ is the groundstate energy of the one-well operator. We conclude by noting that $S_\eta \to S_0$ as $\eta \to 0$.
	\end{proof}

\section{The interaction term}\label{sec.4}
The aim of this section is to prove the following propositions, which follow from the celebrated strategy developped by Helffer and Sjöstrand in \cite{HS84, HS85b} (see also \cite{Hel88}).
\begin{proposition}\label{prop.w}
	We have
	\[\lambda_2(h)-\lambda_1(h)=2|w_h|+\tilde{\mathscr{O}}(e^{-2S_\eta/h})\,,\]
	where $w_h=\langle(\mathscr{L}_h-\mu)\tilde\phi_\ell,\tilde\phi_r\rangle$, $S_\eta$ given in \eqref{eq.Seta}, and $\tilde \phi_\ell$, $\tilde \phi_r$ in \eqref{eq.phitilde}.
\end{proposition}

Hence, the splitting is given by the coefficient $w_h$. Of course, we have to prove that the remainders are smaller than $w_h$, i.e. $2S_0 > S$, which is done in Section~\ref{sec.5}, and requires $L > (2+\sqrt{6})a$. Similarly to the standard strategy, we can rewrite the coefficient $w_h$ as an integral on the line $x_1=0$ separating the two wells. The novelty lies in the emergence of complex coefficients and phases.

\begin{proposition}\label{prop.wh}
Letting $\sigma_\star(x)=\int_{[0,x]}(\mathbf{A}-\mathbf{A}_\star)\cdot\mathrm{d}s$ and 
\begin{equation*}
	\theta(x)=\sigma_r(x)-\sigma_\ell(x)=\int_{[0,x]}(\mathbf{A}_{\ell}-\mathbf{A}_{r})\cdot\mathrm{d}s,
\end{equation*}
we have
	\begin{equation}\label{eq.w0}
	w_h=h^2\int_{\R}\left(2\partial_1\phi_\ell(0,x_2)\,{\phi_\ell}(0,x_2)+i\frac{k(0,x_2)}{h}\phi_\ell(0,x_2)^2\right)e^{-i\theta(0,x_2)/h}\mathrm{d}x_2\,.
	\end{equation}
	with $k(x)=\partial_1\sigma_\ell(x)+\partial_1\sigma_r(x).$
	Moreover, we have explicit formulas
	\begin{equation}\label{eq:theta0}
	\theta(0,x_2)=\frac{b_1L x_2}{2}+2M\arctan(2x_2/L)\,,
	\end{equation}
	and
	\begin{equation}\label{eq:k0}
	k(0,x_2)=b_1 x_2+\frac{2Mx_2}{L^2/4+x_2^2}\,,
	\end{equation}
	where $M$ is given in \eqref{def.M0}.
\end{proposition}
Note that the quantity $M$ can be rewritten in terms of $\beta$, 
\begin{equation}\label{eq.M}
M = \int_0^{\infty} \big( \beta(u) - b_1 \big) \mathrm d u.
\end{equation}
 For the convenience of the reader, we will sum up the strategy in the following sections and try to make it short and transparent.

\subsection{Proof of Proposition \ref{prop.w}}
We let $F=\ker(\mathscr{L}_h-\lambda_1(h))\oplus\ker(\mathscr{L}_h-\lambda_2(h))$ and let $\Pi_F$ be the associated projection.
Let us consider the projections of our one-well quasimodes $g_\star=\Pi_F \tilde\phi_\star$. We let
\[G=\begin{pmatrix}
\langle g_\ell,g_\ell\rangle&\langle g_\ell,g_r\rangle\\
\langle g_r,g_\ell\rangle&\langle g_r,g_r\rangle
\end{pmatrix}=\begin{pmatrix}
g_\ell\\
g_r
\end{pmatrix}\cdot\begin{pmatrix}
g_\ell&g_r
\end{pmatrix}\geq 0\,.\]
We consider the matrix of the quadratic form associated with $\mathscr{L}_h$ in the basis $(g_\ell,g_r)$,
\[L=\begin{pmatrix}
	\langle \mathscr{L}_hg_\ell,g_\ell\rangle&\langle  \mathscr{L}_hg_\ell,g_r\rangle\\
	\langle  \mathscr{L}_hg_r,g_\ell\rangle&\langle  \mathscr{L}_hg_r,g_r\rangle
\end{pmatrix}.\]
The following lemma describes the asymptotic behavior of the matrix $L$, which essentially follows from the classical estimates in \cite{HS84} (see also the pedagogical article \cite{BHR17}).

\begin{lemma}
We have
\begin{equation}\label{eq.approx}
\|g_\star-\tilde\phi_\star\|=\tilde{\mathscr{O}}(e^{-S_\eta/h})\,,\quad Q_h(g_\star-\tilde\phi_\star)=\tilde{\mathscr{O}}(e^{-2S_\eta/h})\,.
\end{equation}
Moreover, the family $(g_r,g_\ell)$ is asymptotically an orthonormal basis of $F$. More precisely, we have
\begin{equation}\label{eq.G}
G=\mathrm{Id}+\mathrm{T}+\tilde{\mathscr{O}}(e^{-2S_\eta/h})\,,\quad \mathrm{T}=\begin{pmatrix}
0&\langle \tilde\phi_\ell,\tilde\phi_r\rangle\\
\langle \tilde\phi_r,\tilde\phi_\ell\rangle&0
\end{pmatrix}=\tilde{\mathscr{O}}(e^{-S_\eta/h})\,.
\end{equation}
In addition, we have
\begin{equation}\label{eq.L}
L=\begin{pmatrix}
	\mu_h& w_h\\
	\overline{w}_h&\mu_h
	\end{pmatrix}+\mu_h\mathrm{T}+\tilde{\mathscr{O}}(e^{-2S_\eta/h})\,,\quad w_h=\langle(\mathscr{L}_h-\mu)\tilde\phi_\ell,\tilde\phi_r\rangle\,.
	\end{equation}

\end{lemma}
\begin{proof}
We recall that
\[(\mathscr{L}_h-\mu_h)\tilde\phi_\star=\tilde{\mathscr{O}}(e^{-S_\eta/h})\,,\]
and notice that, by definition of $F$ and Proposition \ref{prop.goodub},
\[(\mathscr{L}_h-\mu_h)g_\star=\tilde{\mathscr{O}}(e^{-S_0/h})\,.\]	
We infer that
\begin{equation}\label{eq.quasidiff}
(\mathscr{L}_h-\mu_h)(g_\star-\tilde\phi_\star)=\tilde{\mathscr{O}}(e^{-S_\eta/h})\,.
\end{equation}
The first $L^2$-estimate in \eqref{eq.approx} follows by using the spectral theorem and the spectral gap given in Proposition \ref{prop.specgap}. Then, we get the second estimate by using \eqref{eq.quasidiff} and taking the scalar product with $g_\star-\tilde\phi_\star$.

Let us now discuss the asymptotic behavior of $G$. By the Pythagorean theorem, we have
\[\langle \tilde\phi_\star,\tilde\phi_\star\rangle=\langle g_\star,g_\star\rangle+\langle\tilde\phi_\star-g_\star,\tilde\phi_\star-g_\star\rangle\,.\]
By using the support of $\chi_\star$, the definition of $\phi_\star$ and its exponential decay, we have
\begin{equation}\label{eq.tildephi2=1}
\langle \tilde\phi_\star,\tilde\phi_\star\rangle=1+\tilde{\mathscr{O}}(e^{-2S_\eta/h})\,.
\end{equation}
With \eqref{eq.approx}, we deduce that
\[\langle g_\star,g_\star\rangle=1+\tilde{\mathscr{O}}(e^{-2S_\eta/h})\,.\]
In the same way, we write
\[\langle \tilde\phi_\ell,\tilde\phi_r\rangle=\langle g_\ell,g_r\rangle+\langle\tilde\phi_\ell-g_r,\tilde\phi_\ell-g_r\rangle\,,\]
and we notice that
\[\langle \tilde\phi_\ell,\tilde\phi_r\rangle=\tilde{\mathscr{O}}(e^{-S_\eta/h})\,.\]
We deduce \eqref{eq.G}.

Let us finally deal with $L$. The analysis is similar. First, we write
\[\langle \mathscr{L}_h\tilde\phi_\star,\tilde\phi_\star\rangle=\langle \mathscr{L}_hg_\star,g_\star\rangle+\langle\mathscr{L}_h(\tilde\phi_\star-g_\star),\tilde\phi_\star-g_\star\rangle\,.\]
With the localization formula (sometimes called "IMS formula", see \cite[Prop. 4.2]{Raymond}) and using again \eqref{eq.tildephi2=1}, we get
\[\langle \mathscr{L}_h\tilde\phi_\star,\tilde\phi_\star\rangle=\mu_h+\tilde{\mathscr{O}}(e^{-2S_\eta/h})\,,\]
so that, with \eqref{eq.approx},
\[\langle \mathscr{L}_hg_\star,g_\star\rangle=\mu_h+\tilde{\mathscr{O}}(e^{-2S_\eta/h})\,.\]
Then, we write
\[\langle \mathscr{L}_h\tilde\phi_\ell,\tilde\phi_r\rangle=\langle \mathscr{L}_hg_\ell,g_r\rangle+\langle\mathscr{L}_h(\tilde\phi_\ell-g_r),\tilde\phi_\ell-g_r\rangle\,,\]
and we get
\[\langle \mathscr{L}_hg_\ell,g_r\rangle=\langle \mathscr{L}_h\tilde\phi_\ell,\tilde\phi_r\rangle+\tilde{\mathscr{O}}(e^{-2S_\eta/h})\,,\]
where we used \eqref{eq.approx} and the Cauchy-Schwarz inequality to control the remainder. We deduce that
\[\langle \mathscr{L}_hg_\ell,g_r\rangle=\mu_h\langle\tilde\phi_\ell,\tilde\phi_r\rangle+\langle (\mathscr{L}_h-\mu_h)\tilde\phi_\ell,\tilde\phi_r\rangle+\tilde{\mathscr{O}}(e^{-2S_\eta/h})\,,\]
and \eqref{eq.L} follows.
	\end{proof}

We have all the elements at hand to prove Proposition \ref{prop.w}.
The matrix $L$ is the matrix of $\mathscr{L}_h$ in a non-orthonormal basis of $F$. That is why we consider the new family
\[\begin{pmatrix}
\mathfrak{g}_\ell\\
\mathfrak{g}_r
\end{pmatrix}=G^{-\frac12}\begin{pmatrix}
g_\ell\\
g_r
\end{pmatrix}\,,\]
which is orthonormal since
\[\begin{pmatrix}
	\mathfrak{g}_\ell\\
	\mathfrak{g}_r
\end{pmatrix}\cdot\begin{pmatrix}
	\mathfrak{g}_\ell&\mathfrak{g}_r
\end{pmatrix}=G^{-\frac12}\begin{pmatrix}
g_\ell\\
g_r
\end{pmatrix}\cdot\begin{pmatrix}
g_\ell&g_r
\end{pmatrix}G^{-\frac12}=\mathrm{Id}\,.\]
Now, the matrix of $\mathscr{L}_h$ in the basis $(\mathfrak{g}_\ell,\mathfrak{g}_r)$ is $G^{-\frac12}LG^{-\frac12}$ and we have
\[\begin{split}
	G^{-\frac12}LG^{-\frac12}&=\left(1-\frac{\mathrm{T}}{2}\right)\left(\begin{pmatrix}
	\mu_h& w_h\\
	\overline{w}_h&\mu_h
\end{pmatrix}+\mu_h\mathrm{T}\right)\left(1-\frac{\mathrm{T}}{2}\right)+\tilde{\mathscr{O}}(e^{-2S_\eta/h})\\
&=\begin{pmatrix}
	\mu_h& w_h\\
	\overline{w}_h&\mu_h
\end{pmatrix}+\tilde{\mathscr{O}}(e^{-2S_\eta/h})\,.
\end{split}\]
The proposition follows, since the eigenvalues of this matrix are $\mu_h \pm |w_h|$.

\subsection{Proof of Proposition \ref{prop.wh}}

Let us first establish the integral formula \eqref{eq.w0} for the interaction term $w_h$. Recalling the definition \eqref{eq.phitilde} of $\tilde{\phi}_\ell$ and $\tilde{\phi}_r$, we have
\[w_h=\langle(\mathscr{L}_h-\mu_h)\chi_\ell\hat\phi_\ell,\chi_r\hat\phi_r\rangle=\langle[\mathscr{L}_h,\chi_\ell]\hat\phi_\ell,\chi_r\hat\phi_r\rangle\,.\]	
Thus, by letting $P=-ih\nabla-\mathbf{A}$, we get
\[w_h=\langle(P\cdot[P,\chi_\ell]+[P,\chi_\ell]P)\hat\phi_\ell,\chi_r\hat\phi_r\rangle=\langle(P\cdot[P,\chi_\ell]+[P,\chi_\ell]P)\hat\phi_\ell,\hat\phi_r\rangle_{L^2(\R^2_+)}\,.\]
By partial integration we deduce
\[\begin{split}
	w_h&=\langle [P,\chi_\ell]\hat\phi_\ell,P\hat\phi_r\rangle_{L^2(\R^2_+)}+\langle [P,\chi_\ell]P\hat\phi_\ell,\hat\phi_r\rangle_{L^2(\R^2_+)}\\
	&=-ih\int_{x_1>0}\chi'_\ell(x_1)\left(\hat\phi_\ell\overline{(-ih\partial_1-A_1)\hat\phi_r}+(-ih\partial_1-A_1)\hat\phi_\ell\overline{\hat\phi_r}\right)\mathrm{d}x\,.
	\end{split}\]
	Our choice of gauge is such that $A_1=0$. Then, again by partial integration,
	\begin{multline}
	w_h=h^2\int_{\R}\left(\partial_1\hat\phi_\ell\,\overline{\hat\phi_r}-\hat\phi_\ell\,\overline{\partial_1\hat\phi_r}\right)(0,x_2)\mathrm{d}x_2 \\
	+h^2\int_{x_1>0}\chi_\ell\left(\partial_1(\partial_1\hat\phi_\ell\,\overline{\hat\phi_r})-\partial_1(\hat\phi_\ell\,\partial_1\overline{\hat\phi_r})\right) \mathrm{d}x\,. 
		\end{multline}
We notice that
\[h^2\int_{x_1>0}\chi_\ell\left(\partial_1(\partial_1\hat\phi_\ell\,\overline{\hat\phi_r})-\partial_1(\hat\phi_\ell\,\partial_1\overline{\hat\phi_r})\right) \mathrm{d}x=\int_{x_1>0}\chi_\ell\left(\hat\phi_\ell\,\overline{\mathscr{L}_h\hat\phi_r}-\mathscr{L}_h\hat\phi_\ell\,\overline{\hat\phi_r}\right)\mathrm{d}x=0\,,\]
	where we used a partial integration with respect to $x_2$ and the fact that, on the support of $\chi_\ell$, $\mathscr{L}_h\hat\phi_r=\mu_h\hat \phi_r$ and $\mathscr{L}_h\hat\phi_\ell=\mu_h\hat \phi_\ell$. There remains to notice that $\hat\phi_\star=e^{i\sigma_\star/h}\phi_\star$, and with the symmetry $\phi_\ell(x_1,x_2) = \phi_r(-x_1,x_2)$ arround $x_1=0$ we find \eqref{eq.w0}.
	
We end this section by giving explicit formulas for $\theta$ and $k$, which concludes the proof of Proposition \ref{prop.wh}.

\begin{lemma}\label{lem:theta}
	We have
	\begin{equation}\label{eq:theta02}
	\theta(0,x_2)=\frac{b_1L x_2}{2}+2M\arctan(2x_2/L)\,,
	\end{equation}
	and
	\begin{equation}\label{eq:k02}
	k(0,x_2)=b_1 x_2+\frac{2Mx_2}{L^2/4+x_2^2}\,.
	\end{equation}
\end{lemma}

\begin{proof}
	For all $x$ such that $x_1\in(-L/2+a, L/2-a)$,	we have $\nabla\times(\mathbf{A}-\mathbf{A}_\star)=0$, and thus
	\[\begin{split}\sigma_\star(x)&=\int_{[0,x]}(\mathbf{A}-\mathbf{A}_\star)\cdot \mathrm{d}s\\
	&=\int_{0}^{x_1}(\mathbf{A}-\mathbf{A}_\star)_1(u,0)\mathrm{d}u+\int_{0}^{x_2}(\mathbf{A}-\mathbf{A}_\star)_2(x_1,v)\mathrm{d}v\\
	&=\int_{0}^{x_2}(\mathbf{A}-\mathbf{A}_\star)_2(x_1,v)\mathrm{d}v\,,
	\end{split}\]
	where we used the explicit expression \eqref{eq.Al} of $\mathbf{A}_\star$ and the fact that $A_1=0$.
	We have
	\[A_\ell(0,x_2)=\int_0^1\beta\Big(\frac{t^2|x-x_\ell|^2}{2}\Big)t\mathrm{d}t\begin{pmatrix}
		-x_2\\
		\frac{L}{2}
	\end{pmatrix}=\frac{1}{|x-x_\ell|^2}\int_0^{|x-x_\ell|^2/2}\beta(u)\mathrm{d}u\,\begin{pmatrix}
		-x_2\\
		\frac{L}{2}
	\end{pmatrix}\,.\]
	Thus,
	\[A_{\ell,1}(0,x_2) =- \frac{b_1 x_2}{2} - \frac{Mx_2}{L^2/4 + x_2^2}, \qquad A_{\ell,2}(0,x_2)=\frac{b_1L}{4}+\frac{LM/2}{x_2^2+L^2/4}\,.\]
	By symmetry, $\theta$ is given by
	\[\theta(0,x_2) = (\sigma_r - \sigma_\ell)(0,x_2) = 2\int_0^{x_2}A_{\ell,2}(0,v)\mathrm{d}v\,,\]
	and we find \eqref{eq:theta02}.	
	For $k(0,x_2)$, we use
	\[\nabla(\sigma_\ell+\sigma_r)=2\mathbf{A}-(\mathbf{A}_\ell+\mathbf{A}_r)\,,\]
	to get
	\[k(0,x_2) = \partial_1(\sigma_\ell+\sigma_r)(0,x_2)=-(\mathbf{A}_\ell+\mathbf{A}_r)_1(0,x_2) =-2 A_{\ell,1}(0,x_2)\,,\]
	and \eqref{eq:k02} follows.
\end{proof}

\begin{remark}
	Note that the term involving $M$ is here due to the variations of the magnetic field. This term is absent in \cite{FSW22, HK22}.	
\end{remark}

\section{Estimate of the interaction integral}\label{sec.5}
This section is devoted to the proof of the following proposition and lemma, which imply Theorem \ref{thm.main} (see Proposition \ref{prop.w}), with the constant $c$ given by
\begin{equation}\label{eq:c}
c = \frac{2}{\Gamma(\delta_0)} \Big( \frac{b_1L^2}{8} \Big)^{\delta_0} \Big( \frac{1-N}{1+\sqrt{1-N}} \Big)^{2\delta_0} (1-N)^{-1/4}\,.
\end{equation}

\begin{proposition}\label{prop.finalw}
The following estimate holds,
\[
w_h  \underset{h\to 0}{\sim} -\frac{1}{\Gamma(\delta_0)} \Big( \frac{b_1L^2}{8} \Big)^{\delta_0} \Big( \frac{1-N}{1+\sqrt{1-N}} \Big)^{2\delta_0} (1-N)^{-1/4}  h^{1-\delta_0} e^{-S/h},
\]
where
\[S=2\Phi_\ell(0)+I=\int_{0}^{L^2/8}\frac{1}{v}\left(\int_0^{v} \beta(u)\mathrm{d}u\right) \mathrm d v+I\,,\]
with
\[I=\frac{b_1 L^2}{4}\left(\frac{N-1}{2}+\sqrt{1-N}-N\ln(1+\sqrt{1-N})\right)>0\,,\]	
where $M$ and $N$ are given in \eqref{eq.M} and \eqref{eq.N} respectively, and 
\begin{equation}\label{eq:def_delta0}
\delta_0 = \frac{b_1 - b_0}{2b_1}.
\end{equation}
\end{proposition}

In view of Proposition \ref{prop.w}, one needs to check that $2S_0>S$. We recall that $S_0$ is given in \eqref{eq.S1}.
\begin{lemma}\label{lem.comparison}
The quantity $I$ can be rewritten as
\begin{equation}
I = \int_{L^2/8}^{\frac{L^2}{8}(1+\sqrt{1-N})^2} \frac{1}{v}\left(\int_0^{v} \beta(u)\mathrm{d}u\right)\mathrm{d}v  - 2 \int_0^{L^2/8} \beta(u) \mathrm{d} u. 
\end{equation}
Moreover, if $L > (2+\sqrt 6) a$ then $2S_0 > S$.
\end{lemma}

\subsection{Description of $\phi_\ell$}
The aim of this section is to give the explicit expression of $\phi_\ell$ in the interaction zone. This can be done by integrating a classical differential equation known as the Kummer's differential equation, which appeared in \cite{FSW22} and \cite{HK22} but with different parameters. Here $\mu_h$ is of order $h$ and the parameter $M$ is not zero since the magnetic field is not constant.

\begin{lemma}\label{lem.explicit}
For $|x-x_\ell|\geq a$,
\[\phi_\ell(x)=C(h)|x-x_\ell|^{\gamma}\int_0^{+\infty}e^{-\frac{b_1}{4h}(1+2t) |x-x_\ell|^2 } m_h(t)\mathrm{d}t\,,\]	
	where
	\begin{equation}\label{eq.gammadeltamh}
	\gamma=\frac{|M|}{h}\,,\quad \delta=\frac{b_1h-\mu_h}{2hb_1}\,,\quad m_h(t)=t^{\delta-1}(1+t)^{\gamma-\delta}\,,
	\end{equation}
	and
	\begin{equation}\label{eq:DefCh}
	C(h)=\frac{\phi_\ell(0)}{|x_\ell|^{\gamma}}\left(\int_0^{+\infty}e^{-\frac{b_1}{4h}(1+2t) |x_\ell|^2 } m_h(t)\mathrm{d}t\right)^{-1}\,.
	\end{equation}
	Note that $\delta \sim \delta_0$ (with $\delta_0$ as defined in \eqref{eq:def_delta0}) as $h \to 0$.
\end{lemma}

\begin{proof}
Since the function $\phi_\ell$ is radial (with center at $x_\ell$), we can write $\phi_\ell(x)=\phi(|x-x_\ell|)$. Then, $\phi$ solves the radial equation
\[-h^2r^{-1}\partial_r r\partial_r\phi+\frac{\alpha(r)^2}{r^2}\phi=\mu_h\phi\,, \quad \text{with } \alpha(r)=\int_0^{r^2/2}\beta(u)\mathrm{d}u\,.\]
We are interested in the region $r\geq a$ where $\beta=b_1$. There, we have
\[\alpha(r)=M+\frac{b_1r^2}{2}\,,\]
where we recall that $M=\int_0^{a^2/2}(\beta-b_1)\mathrm{d}u$. The equation reads
\[-h^2r^{-1}\partial_r r\partial_r\phi+\left(\frac{M}{r}+\frac{b_1 r}{2}\right)^2\phi=\mu_h\phi\,,\]
so that
\[-h^2r^{-1}\partial_r r\partial_r\phi+\left(\frac{M^2}{r^2}+\frac{b^2_1 r^2}{4}\right)\phi=(\mu_h-b_1M)\phi\,.\]
Then, we let $\phi=r^{\gamma}\psi$ and we have
\[-h^2r^{-1}(r^{-\gamma}\partial_r r^{\gamma}) r(r^{-\gamma}\partial_rr^{\gamma})\psi+\left(\frac{M^2}{r^2}+\frac{b^2_1 r^2}{4}\right)\psi=(\mu_h-b_1M)\psi\,,\]
so that
\[-h^2r^{-1}(\partial_r+\gamma r^{-1})r(\partial_r+\gamma r^{-1})\psi+\left(\frac{M^2}{r^2}+\frac{b^2_1 r^2}{4}\right)\psi=(\mu_h-b_1M)\psi\,.\]
We deduce
\[-h^2r^{-1}(\partial_r+\gamma r^{-1})\psi-h^2(\partial_r+\gamma r^{-1})^2 \psi+\left(\frac{M^2}{r^2}+\frac{b^2_1 r^2}{4}\right)\psi=(\mu_h-b_1M)\psi\,.\]
	Then,
	\[-h^2\partial_r^2\psi-h^2(1+2\gamma)r^{-1}\partial_r\psi+\left(\frac{M^2-h^2\gamma^2}{r^2}+\frac{b^2_1 r^2}{4}\right)\psi=(\mu_h-b_1M)\psi\,.\]
We choose $\gamma$ so that $\gamma^2=M^2h^{-2}$.  Thus,
	\[-h^2\partial_r^2\psi-h^2(1+2\gamma)r^{-1}\partial_r\psi+\frac{b^2_1 r^2}{4}\psi=(\mu_h-b_1M)\psi\,.\]
	We let $\rho=c r^2/2$ and we have $c\partial_\rho=r^{-1}\partial_r$. Note that
	\[\partial_r^2=\partial_r r r^{-1}\partial_r=r^{-1}\partial_r+r^2(r^{-1}\partial_r)^2=c\partial_\rho+2c\rho\partial_\rho^2\,.\]
	
	We let $\psi(r)=\Psi(cr^2/2)$ and we find that
	\[\rho\partial_\rho^2\Psi+(1+\gamma)\partial_\rho\Psi-\frac{b^2_1 \rho}{4h^2c^2}\Psi=\frac{b_1M-\mu_h}{2h^2c}\Psi\,,\]
	and then we let
	\[\Psi(\rho)=F(\rho)e^{-\frac{b_1 }{2hc}\rho}.\]
We get
\[\rho\partial_\rho^2 F+\left(1+\gamma-\rho\frac{b_1}{hc}\right)\partial_\rho F=\left(\frac{b_1M-\mu_h}{2h^2c}+\frac{b_1h(1+\gamma)}{2h^2c} \right)F\,.\]	
This leads to choose $\gamma=|M|h^{-1}$ and $c=b_1h^{-1}$ to get
\[\rho\partial_\rho^2 F+\left(1+\gamma-\rho\right)\partial_\rho F=\frac{b_1h-\mu_h}{2hb_1}F\,,\]
which is the well-known Kummer's equation (and solvable by means of the Fourier or Laplace transform).	Therefore, since we are looking for a decaying solution, we have, for some normalizing constant $C(h)$,
\[F(\rho)=C(h)\int_0^{+\infty}e^{-t\rho} t^{\delta-1}(1+t)^{\gamma-\delta}\mathrm{d}t\,,\quad \delta=\frac{b_1h-\mu_h}{2hb_1}\,.\]
Therefore,
\begin{align*}
	\phi(r)&=r^{\gamma}e^{-b_1r^2/(4h)}F(b_1r^2/(2h))\\
	&=C(h)r^{\gamma}e^{-b_1r^2/(4h)}\int_0^{+\infty}e^{-b_1 r^2t/(2h)} t^{\delta-1}(1+t)^{\gamma-\delta}\mathrm{d}t\,.
\end{align*}
	\end{proof}	
It is possible to estimate the normalization constant $C(h)$ by means of the WKB-approximation of Proposition~\ref{prop.WKB} and a direct estimate of the integral
	\begin{lemma}\label{lem.Ch}
		The normalization constant satisfies
\[C(h)=h^{-\frac12}a_0(0)e^{-\Phi_\ell(0)/h}e^{b_1 L^2/16h}\frac{\Gamma(\delta_0)^{-1}h^{-\delta_0}}{|x_\ell|^{\frac{|M|}{h}}}\left(\frac{b_1L^2}{8}-|M|\right)^{\delta_0}(1+o(1))\,.\]
where $a_0$ and $\Phi_\ell$ were defined in Proposition~\ref{prop.WKB}.
	\end{lemma}
\begin{proof}
	Letting $\varphi(t)=\frac{b_1L^2}{8}t-|M|\ln(1+t)$, we have
		\[\int_0^{+\infty}e^{-\frac{b_1}{4h}(1+2t) |x_\ell|^2 } m_h(t)\mathrm{d}t=e^{-\frac{b_1L^2}{16h}}\int_0^{+\infty}e^{-\varphi(t)/h}t^{\delta-1}(1+t)^{-\delta}\mathrm{d}t\,.\]
Since $\frac{b_1L^2}{8}-|M|>0$ (by \eqref{eq.N}), we have $\varphi'>0$. Note also that $\varphi$ is strictly convex. By means of the change of variable $u=\varphi(t)$, we get
	\[\int_0^{+\infty}e^{-\varphi(t)/h}t^{\delta-1}(1+t)^{-\delta}\mathrm{d}t=(1+o(1))\frac{\Gamma(\delta_0) h^{\delta_0}}{\left(\frac{b_1L^2}{8}-|M|\right)^{\delta_0}}\,.\]
	By inserting in \eqref{eq:DefCh} we get
	\[C(h)=(1+o(1))h^{-\delta_0}\Gamma(\delta_0)^{-1}e^{b_1 L^2/16h}\frac{\phi_\ell(0)}{|x_\ell|^{\frac{|M|}{h}}}\left(\frac{b_1L^2}{8}-|M|\right)^{\delta_0}\,.\]
	The result now follows by inserting the approximation $\phi_\ell(0) \approx \phi^{\mathrm{WKB}}_\ell(0)$, i.e. the WKB-approximation from
	Proposition~\ref{prop.WKB}.
	\end{proof}	

	\subsection{Consequence of the explicit expression}
	Given Lemma~\ref{lem.explicit}, the function $w_h$ is an explicit integral. Let us analyze its asymptotic behavior and establish Proposition \ref{prop.finalw}. First note, with Lemma~\ref{lem.explicit}, that
	\begin{equation*}
	\phi_\ell(0,x_2) = C(h) \int_0^{\infty} \Big( \frac{L^2}{4} + x_2^2 \Big)^{\frac{\gamma}{2}}e^{-\frac{b_1}{4h}(1+2t) (\frac{L^2}{4}+x_2^2)}m_h(t)\mathrm{dt},
	\end{equation*}
	and {\footnotesize
	\begin{equation*}
			\partial_1\phi_\ell(0,x_2)
			=C(h)\int_0^{\infty} \Big(\frac{\gamma L/2}{\frac{L^2}{4}+x_2^2}-\frac{b_1L}{4h}(1+2t)\Big) \Big(\frac{L^2}{4}+x_2^2 \Big)^{\frac{\gamma}{2}}e^{-\frac{b_1}{4h}(1+2t) (\frac{L^2}{4}+x_2^2)}m_h(t)\mathrm{dt}\,.
		\end{equation*}}

\subsubsection{Rescaling of the interaction integral}
		In \eqref{eq.w0} we let $x_2=\frac{Ly}{2}$ and we have
		\[	w_h=\frac{L}{2}h^2\int_{\R}\Big(2\partial_1\phi_\ell(0,\frac{Ly}{2})\,{\phi_\ell}(0,\frac{Ly}{2})+i\frac{k(0,\frac{Ly}{2})}{h}\phi_\ell(0,\frac{Ly}{2})^2\Big) e^{-i\theta(0,\frac{Ly}{2})/h}\mathrm{d}y\,.\]
		Note that Lemma \ref{lem:theta} gives
	\begin{equation*}
			\theta(0,\frac{Ly}{2})=\frac{b_1 L^2}{4}y+2M\arctan y\,, \qquad
	k(0,Ly/2)=\frac{b_1 L}{2}\Big(y-\frac{Ny}{1+y^2}\Big)\,,
		\end{equation*}
	with $N$ from \eqref{eq.N}. Moreover, with $\langle y \rangle^2 = 1 + y^2$,
	\[\phi_\ell(0,Ly/2)=C(h)\Big(\frac{L}{2}\Big)^{\gamma}\langle y\rangle^\gamma\int_0^{+\infty}e^{-\frac{b_1L^2}{16h}(1+2t) \langle y\rangle^2 } m_h(t) \mathrm{d}t\,,\]	
	and	{\footnotesize
		\begin{equation*}
	\partial_1\phi_\ell\Big(0,\frac{Ly}{2}\Big)
	=C(h)\Big(\frac{L}{2}\Big)^\gamma\int_0^{+\infty} \Big(\frac{2\gamma  /L}{1+y^2}-\frac{b_1L}{4h}(1+2t)\Big)\langle y\rangle^\gamma e^{-\frac{b_1 L^2}{16h}(1+2t) \langle y \rangle^2}m_h(t)\mathrm{dt}\,.
\end{equation*}}

Therefore, recalling the definitions of $M, N$ from \eqref{def.M0} and \eqref{eq.N} we can write
\begin{equation}\label{eq.WW}
 w_h=h C(h)^2 \Big(\frac{L}{2}\Big)^{2\gamma} \frac{b_1L^2}{4}W_h\,,
\end{equation}
where (with $s=(s_1,s_2) \in \R^2$)
\begin{equation}\label{eq.WWW}
W_h=\int_{\mathbb{R}}\int_{\R_+\times\R_+}\omega(s,y)\langle y\rangle^{2\gamma}m_h(s_1)m_h(s_2) e^{-\frac{b_1L^2}{8h}(1+s_1+s_2)\langle y\rangle^2 -\frac i h \theta(0,\frac{Ly}{2})}\mathrm{d}s\mathrm{d}y
\end{equation}
and
\begin{equation}\label{eq.omega}
\omega(s,y)=\left(\frac{N}{1+y^2}-1\right)(1-iy)-2s_1\,.
\end{equation}
Now we want to estimate $W_h$.
\subsubsection{The phase}\label{sec.phasey}
In view of \eqref{eq.WWW}, it is natural to introduce the following complex phase
\[\psi(s,y)=\frac{b_1L^2}{8}(1+s_1+s_2)\langle y\rangle^2+i\frac{b_1 L^2}{4}y-2i|M|\arctan y-|M|\ln(1+y^2)\,,\]
where the logarithmic term is extracted from $\langle y \rangle^\gamma$.
This may also be written as
\begin{equation}\label{eq.psi}
\psi(s,y)=\frac{b_1L^2}{8}c(s)\langle y\rangle^2+i\frac{b_1L^2}{4} y-2|M|\ln(1+iy)\,,
\end{equation}
with
\begin{equation}
c(s):= 1 + s_1 + s_2.
\end{equation}

Therefore, \eqref{eq.WWW} becomes
\begin{equation}\label{eq.WWWW}
	W_h=\int_{\mathbb{R}}\int_{\R_+\times\R_+}\omega(s,y)m_h(s_1)m_h(s_2) e^{-\psi(s,y)/h}\mathrm{d}s\mathrm{d}y\,,
\end{equation}
where $\omega$ is given in \eqref{eq.omega} and $m_h$ in \eqref{eq.gammadeltamh}.

\subsubsection{Critical points in $y$}\label{sec.crit}
We compute the complex critical points of $\psi$. 
\begin{lemma}\label{lem.Red2psi>0}
	 Letting 
	\[z_-=z_-(s)=\frac{c(s)-1}{2c(s)}-\frac{1}{2}\sqrt{\left(\frac{c(s)-1}{c(s)}\right)^2+\frac{4}{c(s)}(1-N)}\,,\]	
	we have,  for all $y\in\R$, and all $s\in\R_+\times\R_+$, 
	\[\partial_y\psi(s,iz_-(s))=0\,,\quad \Re \partial^2_y\psi(s,y+iz_-)\geq\frac{b_1L^2}{4}(1-N)\,,\quad \Im \partial^2_y\psi(s,iz_-)=0\,.\]
	In particular, $\Re \psi(s,\cdot+iz_-)$ is strongly convex uniformly in $s$.
	We also have $|z_-(s)|\leq \frac{1}{c(s)} \leq 1$ and $z_-(0)=-\sqrt{1-N}$.
\end{lemma}

\begin{proof}
	The complex critical points of $y \mapsto \psi(s,y)$ solve the equation
	\[c(s)y+i-\frac{iN}{1+iy}=0\,.\]
	With $y=iz$, we get
	\[z^2-\frac{c-1}{c}z+\frac{N-1}{c}=0\,.\]
	We recall that $N\in(0,1)$. The discriminant is equal to
	\[\Delta=\Big(\frac{c-1}{c}\Big)^2+\frac{4}{c}(1-N)\,,\]
	and it is positive. In this case, there are two real roots:
	\[z_\pm=\frac{c-1}{2c}\pm\frac{1}{2}\sqrt{\Big(\frac{c-1}{c}\Big)^2+\frac{4}{c}(1-N)}\,.\]
	We let $y_\pm=iz_\pm$. Since $N \in (0,1)$, we have $z_-\in(- c^{-1}, 0)$.
	It is interesting to notice that
	\[\partial^2_y\psi(s,y)=\frac{b_1 L^2}{4} \Big( c-\frac{N}{(1+iy)^2} \Big)\,,\]
	so that, for all $y\in\R$,
	\[\partial^2_y\psi(s,y+iz_-)=\frac{b_1 L^2}{4} \Big(c(s)-N\frac{(1-z_--iy)^2}{((1-z_-)^2+y^2)^2} \Big)\,,\]
	and thus
	\[\partial^2_y\psi(s,y+iz_-)=\frac{b_1 L^2}{4} \Big( c(s)-N\frac{(1-z_-)^2-y^2}{((1-z_-)^2+y^2)^2}+2iN\frac{(1-z_-)y}{((1-z_-)^2+y^2)^2} \Big)\,.\]
	The imaginary part vanishes at $y=0$, and the real part is bounded from below, which concludes the proof.
	\end {proof}
	We also have the following useful property of the critical point.
	\begin{lemma}\label{lem.dsjz-}
		For all $s\in\R_+\times\R_+$, we have $\partial_{s_j}z_-(s)>0$ and
		\[\partial_{s_j}z_-(0)=\frac{z_-(0)+N-1}{2z_-(0)}= \frac 1 2 (1+\sqrt{1-N})\,.\]	
	\end{lemma}
	\begin{proof}
		We let
		\[P(c,z)=z^2-\frac{c-1}{c}z+\frac{N-1}{c}\,,\]
		and we notice that
		$P(c(s),z_-(s))=0.$	
		In particular, we have
		\[\partial_c P+(\partial_{s_j}z_-)\partial_z P=0\,.\]
		The result follows since $\partial_z P(c(s),z_-(s))=2z_-(s)-\frac{c(s)-1}{c(s)}<0$ and $\partial_c P(c(s),z_-(s))=-(z_-(s)+N-1)c(s)^{-2}>0$.
	\end{proof}
	\subsubsection{Critical points in $s$}\label{sec.crit2}
	Let us look at \eqref{eq.WWWW} and underline that $m_h$ depends on $h$ (see \eqref{eq.gammadeltamh}). This suggests to consider the function
	\begin{equation}\label{eq:def_F}
	F(s)= \psi(s,iz_-(s))-|M|\ln[(1+s_1)(1+s_2)]\,,
	\end{equation}
	which is well-defined and smooth in a neighborhood of $\R_+\times\R_+$.
	\begin{lemma}\label{lem.psi1}
		The function $[0,+\infty)^2\ni s\mapsto F(s)$ is real-valued, it tends to $+\infty $ at infinity and it is strictly convex. It has a unique minimum at $(0,0)$, whose value is
		\[F(0)=\frac{b_1 L^2}{4}\Big(\frac{N}{2}+\sqrt{1-N}-N\ln(1+\sqrt{1-N})\Big) \geq\frac{b_1 L^2}{8}>0\,.\]
		Moreover, $\nabla F(0)=0$ and $\nabla^2F(0)>0$.
	\end{lemma}
	\begin{proof}
		By \eqref{eq.psi} we have
		{\small
		\begin{equation*}
			F(s)=\frac{b_1 L^2}{4} \Big( \frac{ c(s)}{2}(1-z_-^2(s))- z_-(s)-N\ln(1-z_-(s))\Big)
			-|M|\ln[(1+s_1)(1+s_2)] ,
		\end{equation*}
		}
		which is real for all $s\in\R_+\times\R_+$.
		Then, since $iz_-$ is a critical point in $y$, we have  \[\partial_s\left(\psi(s,iz_-(s))\right)=\partial_s\psi(s,iz_-(s))+i\partial_y\psi(s,iz_-(s))\partial_sz_-=(\partial_s\psi)(s,iz_-(s))\,.\]
		Moreover, we have
		\[\partial_s\psi(s,y)=\frac{b_1 L^2}{8}\langle y\rangle^2\begin{pmatrix}
			1\\
			1
		\end{pmatrix}\,.\]
		Hence,
		\[\nabla_s F(s)=\frac{b_1 L^2}{8}(1-z^2_-(s))\begin{pmatrix}
			1\\
			1
		\end{pmatrix}-|M|\begin{pmatrix}
			\frac{1}{1+s_1}\\
			\frac{1}{1+s_2}
		\end{pmatrix}\,.\]
		Note that $z_-(0)=-\sqrt{1-N}$  and thus $\nabla_sF(0)=0$. Computing the second derivatives and using that $\partial_{s_1}z_-=\partial_{s_2}z_->0$, we see that, for all $s$, \[\partial_{s_1}\partial_{s_2}F = \partial^2_{s_j}F  - \frac{| M |}{(1+s_j)^2}>0.\] Thus, $\nabla^2F>0$. This implies that the derivative along any line starting from the origin is increasing, and the result follows.
	\end{proof}
	
	\subsubsection{Estimate of $W_h$}\label{sec.w1w2}	
We recall \eqref{eq.WWWW} and we write
	\[W_h=\int_{\R_+\times\R_+}P(s)e^{-F(s)/h}\int_{\mathbb{R}}e^{-(\psi(s,y)-\psi(s,iz_-(s)))/h}\omega(s,y)\mathrm{d}y\mathrm{d}s\,,\]
	with $P(s)=(s_1s_2)^{\delta-1}[(1+s_1)(1+s_2)]^{-\delta}$ and $F$ defined in \eqref{eq:def_F}.

\begin{lemma}\label{lem:contour}
	For all $s \in \R_+ \times \R_+$ we have
	\begin{equation}\label{eq:contour}
		\int_\R e^{-\psi(s,y)/h} \omega(s,y) \mathrm{d} y = \int_{\R} e^{-\psi(s,y+iz_-(s))/h} \omega(s,y+iz_-(s)) \mathrm{d} y\,.
	\end{equation}
\end{lemma}

\begin{proof}
	For $R>0$ very large, we define $\Gamma_R$, the positively oriented rectangle in the complex plane,
	\[ \Gamma_R = \big[-R,R \big] \cup \big[ -R + i z_- , R + i z_-\big] \cup \big[-R + i z_-, -R \big] \cup \big[R, R + i z_- \big].\]
	Since $| z_-(s) | < 1$, the function $\omega$ is holomorphic inside $\Gamma_R$, and by the Cauchy formula,
	\[ \int_{\Gamma_R}  e^{-\psi(s,y)/h} \omega(s,y) \mathrm{d} y = 0. \]
	The integral on $ \big[-R,R \big]$ (resp. on $ \big[ -R + i z_- , R + i z_-\big]$) converges to the left-hand side (resp. the right-hand side) of \eqref{eq:contour} as $R\to \infty$. Hence, we only have to bound the integral on the two remaining parts, namely $\Gamma_R^{\pm} = \big[\pm R + i z_-, \pm R \big]$. For instance, 
	\[ I^+_R = \int_{\Gamma_R^+}  e^{-\psi(s,y)/h} \omega(s,y) \mathrm{d} y = \int_0^{z_-(s)} e^{-\psi(s,R+it)/h} \omega(s, R+it) i \mathrm{d} t .\]
	For all $t \in [-1,1]$ we can bound
	\[ \mathrm{Re} \, \psi(s,R+it) \geq \frac{b_1 L^2}{4} \big( \frac{R^2}{2} -1 - N \ln \sqrt{4 + R^2} \big), \]
	and we deduce that $I^+_R \to 0$ as $R \to \infty$. We proceed similarly for $\Gamma_R^-$.
\end{proof}
	Considering Lemma \ref{lem:contour} we rewrite
	\[W_h=\int_{\R_+\times\R_+}P(s)e^{-F(s)/h}\int_{\mathbb{R}}e^{-(\psi(s,y+iz_-(s))-\psi(s,iz_-(s)))/h}\omega(s,y+iz_-(s))\mathrm{d}y\mathrm{d}s\,,\]
	and we use the following adaptation of the Laplace method for the $y$-integral.

\begin{lemma}\label{lem:complexlaplace}
For all $s \in \R_+ \times \R_+$ we have
\begin{multline*}
\int_{\R} e^{- (\psi(s,y+iz_-(s)) - \psi(s,iz_-(s)))/h}\omega(s,y+iz_-(s)) \mathrm{d}y \\
= \sqrt{\frac{2\pi h}{\mathrm{Re}\, \partial_y^2 \psi(s,iz_-(s))}} \omega(s,iz_-(s))
+\mathscr{O}(h^{3/2})\,,
\end{multline*}
and the remainder is uniform with respect to $s \in \R_+ \times \R_+$.
\end{lemma}

\begin{proof}
	Let us decompose the phase with its real and imaginary parts:
	\[\psi(s,y+iz_-)-\psi(s,iz_-)=\Re\left(\psi(s,y+iz_-)-\psi(s,iz_-)\right)+i\Theta(s,y)\,,\]
	where $\Theta(s,0)=\partial_y\Theta(s,0)=\partial^2_y\Theta(s,0)=0$ (see Lemma \ref{lem.Red2psi>0}).
	By the Taylor formula, we have
	\[\Re\left(\psi(s,y+iz_-)-\psi(s,iz_-)\right)=y^2\int_0^1(1-u)\Re\partial^2_y\psi(s,iz_-+uy)\mathrm{d}u=\chi_s(y)^2\,,\]
	with
	\[\chi_s(y)=y\Big(\int_0^1(1-u)\Re\partial^2_y\psi(s,iz_-+uy)\mathrm{d}u\Big)^{\frac12}\,,\]
	where we used the positivity of $\Re\partial^2_y\psi(s,iz_-+uy)$ ensured by Lemma \ref{lem.Red2psi>0}. Note that
	\[(\chi_s^{-1})'(0)=\Big(\frac{2}{\Re\partial^2_y\psi(s,iz_-(s))}\Big)^{\frac12}\,.\]
	Thanks to the change of variable $\tilde y=\chi_s (y)/\sqrt{h}$, we find
	\begin{align*}
		\int_{\R} &e^{- (\psi(s,y+iz_-(s)) - \psi(s,iz_-(s)))/h}\omega(s,y+iz_-(s)) \mathrm{d} y \\
		=& \sqrt{h} \int_{\R} e^{-\tilde{y}^2}  \omega(s,iz_-(s)) (\chi_s^{-1})'(0) \mathrm{d} \tilde y\\ &+ \sqrt{h} \int_{\R} e^{- \tilde y^2} (1-e^{-\frac{i}{h} \Theta(s,\chi^{-1}(\sqrt{h}\tilde y))})\omega(s,iz_-(s)) (\chi_s^{-1})'(0) \mathrm{d} \tilde y \\ &+ h^{3/2} \int_{\R} e^{- \tilde y^2} e^{-\frac{i}{h} \Theta(s,\chi^{-1}(\sqrt{h}\tilde y))} \tilde y^2 f_1(\tilde y) \mathrm{d} \tilde y,
	\end{align*}
for some function $f_1$ which is bounded, uniformly with respect to $s$. Here we used that the linear terms in $\tilde y$ give $0$ when integrated against the Gaussian. The first integral gives the main contribution. The third integral is of order $h^{3/2}$. Finally, since $\Theta$ vanishes at the order $3$, we can write $1-e^{-i \Theta/h}$ as a linear term in $\tilde y$ plus a term of order $h$. The linear term vanishes when integrated against the Gaussian, and the second integral is thus of order $h^{3/2}$.
\end{proof}
	
From Lemma \ref{lem:complexlaplace} we see that
	\begin{align}\label{eq.reducedint}
	&\left| W_h - \int_{\R_+^2} \sqrt{\frac{2\pi h}{\mathrm{Re}\, \partial_y^2 \psi(s,iz_-(s))}} P(s)e^{-F(s)/h} \omega(s,iz_-(s)) \mathrm{d}s \right| \nonumber \\
	&\leq C h^{3/2}  \int P(s) e^{-F(s)/h}  \mathrm{d}s 
	\end{align}
We can use the Laplace method to estimate these integrals. We recall from Lemma~\ref{lem.psi1} that $F$ has a unique critical point and global minimum at $s=0$. With \eqref{eq.omega}, we have
\[\omega(s,iz_-(s))=\Big(\frac{N}{1-z_-(s)^2}-1\Big)(1+z_-(s))-2s_1\,.\]
In particular, since $z_-(0)=-\sqrt{1-N}$, 
\[\omega(0,iz_-(0))=0\,.\]	
Thus, we need to consider the linear approximation of $\omega(s,iz_-(s))$ near $s=0$. We have
\begin{align*}
	\omega(s,iz_-(s))&=\Big( \frac{N}{(1-z_-(0))^2}- 1 \Big)\nabla z_-(0)\cdot s-2s_1+\mathscr{O}(|s|^2) \\
	 &=\Big( \frac{N}{(1+\sqrt{1-N})^2} -1 \Big) \frac 1 2 (1 + \sqrt{1-N}) (s_1+s_2) - 2s_1 + \mathscr O (|s|^2) \\
	&= -\sqrt{1-N} (s_1 + s_2) - 2 s_1 + \mathscr O(|s|^2),
\end{align*}
where we used Lemma \ref{lem.dsjz-}.
With the Laplace method, we deduce that
\begin{equation} \label{est.Wh-old}
\int_{\R_+^2} \sqrt{\frac{2\pi h}{\mathrm{Re}\, \partial_y^2 \psi(s,iz_-(s))}} P(s)e^{-F(s)/h} \omega(s,iz_-(s)) \mathrm{d}s=c_0h^{1+\delta_0}e^{-F(0)/h}(1+\mathscr O(\sqrt h))\,,
\end{equation}
with
\begin{equation*}
c_0=-\sqrt{\frac{2\pi}{\Re\partial^2_y\psi(0,iz_-(0))}}
\int_{\R_+^2}(s_1s_2)^{\delta_0-1}\left(\sqrt{1-N} (s_1 + s_2) + 2 s_1\right)e^{-\frac12\mathrm{Hess}_0 F(s,s)}\mathrm{d}s.
\end{equation*}
Note that this integral can be computed explicitly (first replace $2s_1$ by $s_1+s_2$ by symmetry, then on the set $\lbrace s_2>s_1 \rbrace$ use the variables $v=s_2-s_1$ and $u=s_1s_2$). We find
\begin{equation}\label{eq.c0}
c_0 = -\frac{8\pi \Gamma(\delta_0)}{b_1L^2}(1-N)^{-1/4} \Big( \frac{b_1L^2}{8} \Big)^{-\delta_0} (1+ \sqrt{1-N})^{-2\delta_0}.
\end{equation}
The second integral in \eqref{eq.reducedint} can also be calculated using the Laplace method. However, it is easier since the factor $\omega(s,iz_-(s))$ is not present. Since $\omega(s,iz_-(s))$ vanishes linearly at $s=0$ this means that this second integral becomes a $1/\sqrt{h}$ bigger compared to the first, but since it had an extra factor $h$, we get the final result
\begin{equation} \label{est.Wh}
W_h=c_0h^{1+\delta_0}e^{-F(0)/h}\big(1+\mathscr O(\sqrt h)\big)\,,
\end{equation}
with $c_0$ given by \eqref{eq.c0}.

\subsubsection{Estimate of $w_h$: End of the proof of Proposition \ref{prop.wh}}
Recalling \eqref{eq.WW}, we get
\begin{equation*}
w_h = C(h)^2 \Big(\frac{L}{2}\Big)^{2\gamma} \frac{b_1L^2}{4} c_0 h^{2+\delta_0} e^{-F(0)/h}(1+\mathscr O(h^{1/2}))\,.
\end{equation*}
Therefore, with Lemma \ref{lem.Ch},
\begin{equation*}w_h\underset{h\to 0}{\sim}
\frac{b_1L^2}{4} a_0(0)^2e^{-2\Phi_\ell(0)/h} h^{1-\delta_0}\Gamma(\delta_0)^{-2}\Big(\frac{b_1L^2}{8}-|M|\Big)^{2\delta_0}c_0e^{-I/h}\,,
\end{equation*}
where $I=F(0)-\frac{b_1L^2}{8}$. With the explicit value of $F(0)$ in Lemma \ref{lem.psi1} we see that $I>0$. This finishes the proof of Proposition \ref{prop.finalw}.

\subsubsection{Comparison with the remainders and proof of Lemma \ref{lem.comparison}}

We have now to ensure that $S < 2S_0$, i.e.
\begin{equation}\label{eq.inegalite}
\int_{0}^{L^2/8}\frac{1}{v}\left(\int_0^{v} \beta(u)\mathrm{d}u\right)\mathrm{d}v+I<\int_{0}^{(L-a)^2/2}\frac{1}{v}\left(\int_0^{v} \beta(u)\mathrm{d}u\right)\mathrm{d}v\,.
\end{equation}
We can rewrite $I$ using the Agmon distance,
\begin{equation}
I = \int_{L^2/8}^{\frac{L^2}{8}(1+\sqrt{1-N})^2} \frac{1}{v}\left(\int_0^{v} \beta(u)\mathrm{d}u\right)\mathrm{d}v  - 2 \int_0^{L^2/8} \beta(u) \mathrm{d} u. 
\end{equation}
Therefore, \eqref{eq.inegalite} becomes
\begin{equation} \label{eq.inegalite2}
\int_{(L-a)^2/2}^{\frac{L^2}{8}(1+\sqrt{1-N})^2} \frac{1}{v}\left(\int_0^{v} \beta(u)\mathrm{d}u\right)\mathrm{d}v < 2 \int_0^{L^2/8} \beta(u) \mathrm{d} u.
\end{equation}
Note that the left-hand side is positive because
\begin{equation}
\frac L 2 (1 + \sqrt{1-N} ) \geq \frac L 2 (2 - \sqrt{N}) \geq L - a\,,
\end{equation}
where we used \eqref{eq.N}.

Let us assume $L > ( 2 + \sqrt 6 )a$, and prove \eqref{eq.inegalite2} in this case. We bound the left-hand side as follows,
\begin{equation}
\int_{(L-a)^2/2}^{\frac{L^2}{8}(1+\sqrt{1-N})^2} \frac{1}{v}\left(\int_0^{v} \beta(u)\mathrm{d}u\right)\mathrm{d}v \leq \int_{(L-a)^2/2}^{L^2/2} \frac{1}{v} \left(\int_0^{v} \beta(u)\mathrm{d}u \right) \mathrm{d} v.
\end{equation}
and using $\beta \leq b_1$ we deduce
\begin{equation}
\int_{(L-a)^2/2}^{\frac{L^2}{8}(1+\sqrt{1-N})^2} \frac{1}{v}\left(\int_0^{v} \beta(u)\mathrm{d}u\right)\mathrm{d}v \leq \frac{b_1}{2} \big( L^2 - (L-a)^2 \big) = \frac{b_1}{2}\big( 2La - a^2 \big).
\end{equation}
Moreover, the right-hand side of \eqref{eq.inegalite2} is bounded from below by
\begin{equation}
2 \int_0^{L^2/8} \beta(u) \mathrm d u = \frac{b_1L^2}{4} - 2|M| \geq \frac{b_1}{2}\big( \frac{L^2}{2} - 2a^2 \big).
\end{equation}
Thus, \eqref{eq.inegalite2} is true as soon as $L^2/2 - 2a^2 > 2La - a^2$, which is equivalent to $L > (2+\sqrt{6})a$.

\section*{Acknowledgments}
This work was started while N.R. was visiting the University of Copenhagen. He wishes to thank the University of Copenhagen and its QMATH group for their hospitality. This stay was partially funded by the IRN MaDeF (CNRS).
S.F. and L.M. were partially supported by the grant 0135-00166B from the Independent Research Fund Denmark and by the ERC Advanced Grant MathBEC - 101095820.

\bibliographystyle{abbrv}
\bibliography{biblioFMR}
\end{document}